\documentclass[lettersize,journal]{IEEEtran}
\usepackage{amsmath,amsfonts}
\usepackage{algorithmic}
\usepackage{algorithm}
\usepackage{array}
\usepackage[caption=false,font=normalsize,labelfont=sf,textfont=sf]{subfig}
\usepackage{textcomp}
\usepackage{stfloats}
\usepackage{url}
\usepackage{verbatim}
\usepackage{graphicx}
\usepackage{cite}
\usepackage{tikz}
\usepackage{booktabs}

\usepackage[
colorlinks=true,
citecolor=blue,
linkcolor=magenta,
urlcolor=cyan]{hyperref}

\usepackage{listings}

\definecolor{codegreen}{rgb}{0,0.6,0}
\definecolor{codegray}{rgb}{0.5,0.5,0.5}
\definecolor{codepurple}{rgb}{0.58,0,0.82}
\definecolor{backcolour}{rgb}{0.95,0.95,0.92}

\lstdefinestyle{mystyle}{
    backgroundcolor=\color{backcolour},   
    commentstyle=\color{codegreen},
    keywordstyle=\color{magenta},
    numberstyle=\tiny\color{codegray},
    stringstyle=\color{codepurple},
    basicstyle=\ttfamily\small,
    breakatwhitespace=false,         
    breaklines=true,                 
    captionpos=b,                    
    keepspaces=true,                 
    numbers=left,                    
    numbersep=5pt,                  
    showspaces=false,                
    showstringspaces=false,
    showtabs=false,                  
    tabsize=2
}

\definecolor{seagreen}{rgb}{0.18, 0.55, 0.34}
\definecolor{mediumviolet-red}{rgb}{0.78, 0.08, 0.52}
\definecolor{khaki}{rgb}{0.94, 0.9, 0.55}

\lstdefinelanguage{mypython}
{
	keywords=[1]{from, import, assert, not, print},
	keywordstyle=[1]{\color{mediumviolet-red}},
	keywords=[2]{surecr, torch, cp, lo, pl},
	keywordstyle=[2]{\color{seagreen}},
	numbers=none,
	upquote=true,
	showstringspaces=false,
	basicstyle=\ttfamily\tiny,
	columns=fullflexible,
	keepspaces=true,
	emph={True,False,as,def,return,float,class,match,switch,len},
	emphstyle={\color{seagreen}},
	frame=trBL,
	captionpos=b
}

\newcommand{\eg}{{\it e.g.}}
\newcommand{\ie}{{\it i.e.}}

\newcommand{\BEQ}{\begin{equation}}
\newcommand{\EEQ}{\end{equation}}
\newcommand{\BEAS}{\begin{eqnarray*}}
\newcommand{\EEAS}{\end{eqnarray*}}
\newcommand{\ones}{\mathbf 1}
\newcommand{\reals}{{\mbox{\bf R}}}

\newcommand{\Div}{\nabla \!\cdot\!}

\newcommand{\Tr}{\mathop{\bf Tr}}
\newcommand{\diag}{\mathop{\bf diag}}

\newcommand{\Expect}{\mathop{\bf E{}}}

\newcommand{\argmin}{\mathop{\rm argmin}}

\newcommand{\var}{\mathop{\bf var}}

\newcommand{\Card}{\mathop{\bf card}}

\newcommand{\prox}{\mathbf{prox}}

\newcommand{\SURE}{\textrm{SURE}}

\usepackage{relsize} 

\newcommand{\hutchpp}{\text{Hutch\raisebox{0.35ex}{\relscale{0.75}++}}}

\begin{document}

\title{Tractable Evaluation of Stein's Unbiased Risk Estimate with Convex Regularizers}

\author{Parth Nobel,
Emmanuel Cand\`es,~\IEEEmembership{Fellow,~IEEE,}
Stephen Boyd,~\IEEEmembership{Fellow,~IEEE}
\thanks{
Parth Nobel was supported in part by the National Science Foundation
Graduate Research Fellowship Program under Grant No. DGE-1656518. Any
opinions, findings, and conclusions or recommendations expressed in
this material are those of the author(s) and do not necessarily
reflect the views of the National Science Foundation.  

Stephen Boyd was partially supported by ACCESS (AI Chip Center for Emerging
Smart Systems), sponsored by InnoHK funding, Hong Kong SAR, and by Office of
Naval Research grant N00014-22-1-2121.
Emmanuel Cand\`es was supported by the Office of Naval
Research grant N00014-20-1-2157, the National Science Foundation grant
DMS-2032014, the Simons Foundation under award 814641, and the ARO
grant 2003514594.
    }
\thanks{Parth Nobel and Stephen Boyd are with the Department of Electrical Engineering at Stanford University.}
\thanks{Emmanuel Cand\`es is with the Department of Statistics at Stanford University.}
}

\markboth{IEEE Transactions on Signal Processing}
{Nobel \MakeLowercase{\textit{et al.}}: Tractable Evaluation of Stein's Unbiased Risk Estimate with Convex Regularizers}


\maketitle

\begin{abstract}
Stein's unbiased risk estimate (SURE) gives an unbiased estimate of
  the $\ell_2$ risk of any estimator of the mean of a Gaussian random
  vector.  We focus here on the case when the estimator minimizes a
  quadratic loss term plus a convex regularizer.  For these estimators
  SURE can be evaluated analytically for a few special cases, and
  generically using recently developed general purpose methods for
  differentiating through convex optimization problems; these generic
  methods however do not scale to large problems.  In this paper we
  describe methods for evaluating SURE that handle a wide class of
  estimators, and also scale to large problem sizes.
\end{abstract}


\section{Introduction and background}

\subsection{Stein's unbiased risk estimate (SURE)}
We consider $y \sim \mathcal N(\mu, \sigma^2 I)$ where $\mu \in \reals^d$ and
$I$ is the $d \times d$ identity matrix.
We assume $\sigma$ is known and that we are estimating $\mu$.
We are analyzing estimators $\hat\mu: \reals^d \to \reals^d$ which estimate
$\mu$ given a single sample $y$.
The $\ell_2$ risk of an estimator $\hat\mu$ is
$R(\hat\mu) = \Expect{\|\hat\mu(y) - \mu\|_2^2}$.

In 1981, Charles Stein introduced in \cite{stein1981estimation} what
is now called Stein's unbiased risk estimate, \BEQ\label{e-sure}
\SURE(\hat\mu, y) = -d\sigma^2 + \|\hat\mu(y) - y\|_2^2 + 2\sigma^2
\Div \hat\mu(y), \EEQ where
$\Div \hat\mu(y) = \sum_{i=1}^d \frac{\partial \hat\mu_i}{\partial
  y_i}(y)$ is the divergence of $\hat\mu$ at $y$.  The divergence can
also be expressed as $\Div \hat \mu (y) = \Tr(D\hat\mu(y))$, where
$D\hat \mu(y)$ is the $d\times d$ Jacobian or derivative, evaluated at
$y$, and $\Tr$ denotes the trace of a matrix.  Stein showed that the
SURE statistic is an unbiased estimate of the risk in the sense that
$\Expect \SURE(\hat\mu, y) = R(\hat\mu)$.
The challenge in
evaluating $\SURE(\hat \mu,y)$ is evaluating the divergence
$\Div \hat\mu(y)$.

In \eqref{e-sure}, it is assumed that the estimator $\hat \mu$ is weakly 
differentiable and satisfies some integrability conditions.  
If this is not the
case, SURE is not defined;
we discuss this in more detail in \S\ref{s-weak-diff}.

\subsection{Convex regularized regression}
In this paper we consider the setting where $\mu$ is a known linear
function of unknown parameters $\beta \in B$, where $\beta$ can be a
vector, a matrix, or tuples of vectors and matrices, and $B$ is the
vector space of all such parameters, with dimension $p$.  We will
identify $B$ with $\reals^p$, using some fixed ordering of the entries
of the vectors and matrices that comprise $b \in B$.  For $b \in B$,
we define $\|b\|_2^2$ as the sum of the squares of the entries of $b$.
In other words, we use $\|b\|_2^2$ to mean the square of the $\ell_2$
norm of $b$, interpreted as an element of $\reals^p$.  For example, if
$b$ is a matrix, $\|b\|_2^2$ denotes its Frobenius norm, and not its
induced $\ell_2$ norm/maximum singular value.  When $b$ is a matrix
and we wish to refer to its induced $\ell_2$ norm, we use the notation
$\sigma_\text{max}(b)$.

We take $\mu = \mathcal A \beta$, where $\mathcal A: B \to \reals^d$
is linear.
Using our identification  of $B$ and $\reals^p$,
we can represent $\mathcal A$ explicitly as a $d \times p$ matrix.
But for purposes of computing, it is more convenient to keep it abstract.
In the sequel we will denote the adjoint of the mapping as $\mathcal A^*$.

We consider estimators given by convex regularized regression,
\ie, of the form
\BEQ\label{e-estimators}
    \hat\mu(y) = \mathcal A \argmin_b
        \left(\frac{1}{2} \|\mathcal A b - y\|_2^2 + r(b)\right),
\EEQ
where $r: B \to \reals \cup \{\infty\}$ is a convex regularizer.
The data in this problem are the linear mapping $\mathcal A$,
the regularizer $r$, and the observed sample $y$.
We will denote the argmin in \eqref{e-estimators} as $\hat\beta(y)$ so that
$\hat\mu(y) = \mathcal A \hat\beta(y)$.
Many common estimators have this form.
For some of these, there are closed form expressions for either $\hat \mu(y)$
or SURE.

\subsection{This paper}
In this paper we introduce an algorithm to tractably compute SURE for convex 
regularized regression.
Our algorithm, which we call SURE-CR, requires no direct access to the regularizer,
only the ability to evaluate and differentiate its proximal operator, \ie,
a proximal operator oracle.
SURE-CR requires no knowledge of $\mathcal A$ beyond the ability to evaluate
$b \mapsto\mathcal A b$ and $v \mapsto \mathcal A^* v$, \ie, a forward-adjoint
oracle for $\mathcal A$ and $\mathcal A^*$.
SURE-CR easily scales to problems with 
numbers of parameters in the millions,
where forming or storing the matrix $D\hat \mu(y)$ would be impossible.

\subsection{Classical examples of convex regularized regression}\label{s-classical-estimators}
\paragraph{Ordinary least squares}
In ordinary least squares, $\mathcal A$ is a full-rank data matrix $X \in \reals^{d \times p}$
and 
\[
\hat\mu(y) = X \argmin_b \frac{1}{2}\|Xb - y\|_2^2 = X (X^* X)^{-1} X^* y.
\]
With the orthogonal projection matrix $H$ defined as $H = X (X^* X)^{-1} X^*$, we have 
\[
    \SURE (\hat\mu, y) = (2p - d) \sigma^2 + \|H y - y\|_2^2.
\]

\paragraph{Ridge regression}
In ridge regression, $\mathcal A$ is a potentially rank-deficient data matrix, and 
\begin{align*}
    \hat\mu(y) &= X \argmin_b \left( \frac{1}{2}\|Xb - y\|_2^2 +
 \lambda \|b\|_2^2\right) \\
    &= X (X^* X+\lambda I)^{-1} X^* y,
\end{align*}
where $\lambda >0$.
With $H = X (X^* X + \lambda I)^{-1} X^*$, we have
\[
\SURE(\hat\mu,y) = -d \sigma^2 + \|H y - y\|_2^2 + 2\sigma^2\Tr H.
\]

\paragraph{LASSO}
In LASSO, $\mathcal A$ is again a data matrix, and 
\[
\hat\mu(y) = X \argmin_b \left( \frac{1}{2}\|Xb - y\|_2^2 +
\lambda \|b\|_1 \right),
\]
where $\lambda >0$.
There is no analytical formula for $\hat \mu(y)$, but it is readily 
evaluated numerically. In the usual case where the LASSO solution
is unique, SURE takes the form 
\[
\SURE(\hat\mu, y) = -d \sigma^2
    + \|X\hat\beta(y) - y\|_2^2 + 2 \sigma^2 \Card\hat\beta(y),
\]
where $\Card(\cdot)$ is the number of nonzero entries \cite{zou2007degrees}.

The function $\hat\mu$ is non-differentiable on a set of Lebesgue measure $0$.
Therefore, the random data $y$ is almost surely at a differentiable point of $\hat\mu$.
Specifically, one consequence of \cite[Lemma 3]{tibshirani2012degrees} is that $\hat\mu$
is non-differentiable only on the set $\bigcup_{i=1}^p\{z: |(X^T z)_i| = \lambda \}$.

\subsection{Matrix estimators}\label{s-matrix-estimators}
We now describe a few examples where $\mathcal A$ is not a data matrix,
and except for the first example,
there are no known expressions for SURE.

\paragraph{Singular value thresholding}
The first example is
singular value thresholding, where $y$ and $\beta$ are matrices in 
$B= \reals^{m \times n}$ and 
\[
\hat\mu(y) = \argmin_b \left( \frac{1}{2}\|b - y\|_F^2 + 
\lambda \|b\|_* \right),
\]
where $\lambda>0$ and 
$\|\cdot\|_*$ is the nuclear norm, \ie, the dual of the spectral
norm, the sum of the singular values of $b$.  Here we take
$\mathcal A$ to be the identity operator in our generic formulation.
The estimator $\hat \mu$ can be expressed analytically as singular
value thresholding, \ie, $\hat \mu (y) = U F(\Sigma) V^*$,
where
$y=U\Sigma V^*$ is the singular value decomposition of $y$, and
$F(\Sigma)$ is the diagonal matrix with
$F(\Sigma)_{ii} = \max\{\Sigma_{ii}- \lambda,0\}$.  A closed form
expression for SURE in this case is given in \cite{candes2013svt}.

\paragraph{Matrix completion}
Our second example is matrix completion, which extends singular value
thresholding to the setting where only some entries of a matrix are
observed.  As in singular value thresholding, we have
$\beta\in B=\reals^{m \times n}$.  In matrix completion,
$\mathcal A: B \to \reals^d$ is a selection operator, with $d$ the
number of entries of $\beta$ that are being observed (hence, the
observation $\mu$ is a vector containing the observed entries of the
matrix).
(A selection operator is one where each entry of $\mathcal A b$ is
an entry of $b$.)
The estimator is
\[
\hat\mu(y) = \mathcal A\argmin_b \left( \frac{1}{2}\|\mathcal A b
    - y\|_2^2 + \lambda \|b\|_* \right),
\]
where $\lambda>0$.
Unlike singular value thresholding, there is no known analytical
expression for $\hat \mu(y)$, but it is readily evaluated.  Also,
there is no known closed-form expression for SURE for matrix
completion which can be tractably evaluated.

For future use we note that $\hat \beta(y) = 0$ if and only if
\BEQ\label{e-lambda-max}
\lambda \geq \lambda_\text{max} = \sigma_{\max}(\mathcal A^* y), 
\EEQ
where $\mathcal A^* y$ is a matrix which satisfies $\mathcal A \mathcal A^* y = y$
and which has all entries not uniquely determined by that equation equal to $0$.

\paragraph{Robust PCA}
Our final example is robust PCA, where
\[
b = (L, S) \in \reals^{m \times n} \times \reals^{m \times n}
\]
and $\mathcal A (L, S) = L + S$.
For completeness, we note that $\mathcal A^* V = (V, V)$ where $V$ is any matrix.
The estimator is given by 
\[
\hat\mu(y) = \mathcal A \argmin_{L,\;S} \left( \frac{1}{2}
\|\mathcal A (L,S) - y\|_F^2 + \lambda \|L\|_* + \gamma \|S\|_1 \right),
\]
where $\lambda>0$ and $\gamma>0$.
There is no known closed-form expression for $\hat \mu(y)$,
but it is readily evaluated.
There is no known closed-form expression for SURE.

Here too we can determine the values of $\lambda$ and $\gamma$ for which the optimal solution obeys 
$\hat \beta(y) = 0$.
We have $\hat \beta(y) = 0$ if and only if
\BEQ\label{e-rob-pca-zero-soln}
\lambda \geq \lambda_\text{max} = \sigma_\text{max}(y) \quad \mbox{and} \quad 
\gamma \geq \gamma_\text{max} = \|y\|_\infty,
\EEQ
where $\|y\|_\infty = \max_{i,j} | y_{ij}|$.
We are not aware of this result appearing in
the literature, so we give a short derivation here.
The necessary and sufficient optimality condition for $L$ and $S$ is
\[
L+S-y + \lambda \partial \|L\|_* \ni 0,
\qquad
L+S-y + \gamma \partial \|S\|_1 \ni 0,
\]
where $\partial$ denotes the subdifferential.
Applying this to $L=S=0$ we have that $L=S=0$ is optimal if and only if
\[
y \in \lambda \partial \|0\|_*,
\qquad
y \in \gamma \partial \|0\|_1.
\]
Using the fact that the 
subdifferential of a norm at zero is the unit ball of the 
dual norm, we can write this as \eqref{e-rob-pca-zero-soln}.

\subsection{Algorithms for convex regularized regression}\label{s-solver-algorithms}
Several algorithms for evaluating the estimator \eqref{e-estimators}
access the data $\mathcal A$ and $r$ in the following restricted way: the linear operator $\mathcal A$ is accessed only through its forward and adjoint 
oracle. 
This means we can evaluate $\mathcal A b$ for any $b \in B$, and 
$\mathcal A^*z$ for any $z \in \reals^d$, where $\mathcal A^*: \reals^d \to B$ 
is the adjoint of $\mathcal A$.
This allows us to handle problems without forming or storing 
an explicit matrix representation of $\mathcal A$.

The regularizer is accessed only via its proximal operator
$\prox_{\eta r}: B \to B$, given by
\[
    \prox_{\eta r}(v) = \argmin_b \left(\eta r(b)
                                       + \frac{1}{2} \|b - v\|_2^2\right),
\]
where $v, b\in B$ and $t$ is a positive scalar that can be interpreted 
(in the context of algorithms) as a step length.
Thus our access to the regularizer is via
its proximal operator, \ie, we can evaluate $\prox_{\eta r}(v)$ for any $v$.
The proximal operators of many common regularizers are known and readily computed
\cite{beck2017first,chierchia2016proximity,moreau1962fonctions,parikh2014proximal}.

As examples, in LASSO, $r(b) = \lambda \|b\|_1$, and its proximal operator 
is given elementwise by
\[
    (\prox_{\eta r}(v))_i = \begin{cases}
        v_i - \eta \lambda  & \text{if } v_i > \eta \lambda \\
        -v_i + \eta \lambda & \text{if } v_i < \eta \lambda \\
        0 & \text{else}
    \end{cases}.
\]
This function is known as soft-thresholding and we denote it $\mathcal T_{\eta\lambda}$.

In matrix completion, $r(b) = \lambda \|b\|_*$ and $\prox_{\eta r}(v)$ is given by
singular value thresholding with regularization parameter $t \lambda$.
In robust PCA, $r((L, S)) = \lambda \|L\|_* + \gamma \|S\|_1$ is separable with
respect to $L$ and $S$.
Therefore,
\begin{multline*}
        \prox_{\eta r}((L, S))
        = \argmin_{L', S'} \bigg(
        \eta(\lambda \|L'\|_* + \gamma \|S'\|_1) \\
           + \frac{1}{2} \|(L', S') - (L, S) \|_2^2
        \bigg)
    \end{multline*}
    \vspace{-4ex}
    \begin{multline*}
            \phantom{\prox_{\eta r}((L, S))}
        = \argmin_{L', S'} \bigg(
       \eta  \lambda \|L'\|_* + \eta \gamma \|S'\|_1 \\
            + \frac{1}{2} \|L' - L\|_2^2  + \frac{1}{2}\|S' - S\|_2^2
        \bigg)
    \end{multline*}
    \vspace{-4ex}
    \begin{multline*}
            \phantom{\prox_{\eta r}((L, S))}
        = \Bigg(
        \argmin_{L'} \left(\eta \lambda \|L'\|_*+ \frac{1}{2} \|L' - L\|_2^2\right),  \\
        \argmin_{S'} \left(\eta \gamma \|S'\|_1 + \frac{1}{2}\|S' - S\|_2^2\right) 
        \Bigg)
    \end{multline*}
    \vspace{-4ex}
        \begin{multline*}
            \phantom{\prox_{\eta r}((L, S))}
        = \left(\prox_{\eta \,\lambda \|\cdot\|_*}(L),
                 \prox_{\eta \,\gamma\|\cdot\|_1}(S) \right).
        \end{multline*}
These two proximal operators are exactly those in LASSO and matrix completion.

We now mention three algorithms that only require oracle access to $\mathcal A$,
$\mathcal A^*$, and $\prox_{\eta r}(\cdot)$. 

\paragraph{ISTA}
The proximal gradient method (also known as ISTA) 
\cite{beck2017first,beck2009fast,parikh2014proximal}
consists of the iterations
\[
        b^{k+1} = \prox_{\eta r}\left(b^{k}
                        - \eta \mathcal A^* \left( \mathcal A b^{k}
                        -y\right)\right).
\]
The algorithm itself requires only multiplication by $\mathcal A$ and $\mathcal A^*$.
The step length $\eta$ must satisfy
$\eta \leq 2/\sigma_{\max}(\mathcal A)$ to guarantee convergence \cite[\S4.2]{parikh2014proximal}; here, 
$\sigma_{\max}(\mathcal A)$ is the induced $\ell_2$ norm,
which can be computed by a power algorithm 
that only uses multiplication by $\mathcal A$ and $\mathcal A^*$.
For our purposes, ISTA can be initialized with any vector which is selected
independently of $y$,
and we shall require that the Jacobian of $b^1$ with respect to $y$ always be the zero matrix.

\paragraph{FISTA}
The accelerated proximal gradient method (also known as FISTA) 
\cite{beck2017first,beck2009fast,parikh2014proximal}
consists of adding a momentum term to the proximal gradient method to
obtain the iterations
\[
    \begin{aligned}
        \tau^{k+1} &= \frac{1 + \sqrt{1 + 4 \left(\tau^k\right)^2}}{2} \\
        b^{k+1/2} &= b^k + \frac{\tau^k - 1}{\tau^{k+1}} \left(
            b^k - b^{k-1}\right)  \\
        b^{k+1} &= \prox_{\eta r}\left(b^{k+1/2}
                        - \eta \mathcal A^* \left( \mathcal A b^{k+1/2}
                        -y\right)\right),
    \end{aligned}
\]
where $k$ is the iteration counter and $\tau_1 = 1$.
The algorithm itself requires only multiplication by $\mathcal A$ and $\mathcal A^*$.
The step length $\eta$ must satisfy
$\eta \leq 1/\sigma_{\max}(\mathcal A)$ to guarantee convergence \cite[\S4.3]{parikh2014proximal}.
It is also possible to use $\tau^{k} = \frac{k+2}{2}$
\cite[Remark 10.35]{beck2017first}, as we do in the sequel.
For our purposes, FISTA can be initialized with any vector $b^1$ which is selected
independent of $y$,
\ie, we require that the Jacobian of $b^1$ with respect to $y$ always be the zero matrix.
FISTA is almost always preferable to ISTA.

\paragraph{ADMM}
The third algorithm we mention is the alternating direction method of 
multipliers (ADMM) \cite{boyd2011distributed}, 
with iterations 
\[
    \begin{aligned}
    b^{k+1} &= \prox_{\eta r}(z^k - u^k) \\
    z^{k+1} &= (\eta \mathcal A^* \mathcal A + I)^{-1}
                    (b^{k+1} + u^k + \eta \mathcal A^* y)\\
    u^{k+1} &= u^k + b^{k+1} - z^{k+1},
    \end{aligned}
\]
where $u^k,z^k \in B$.
For ADMM, the parameter $\eta$ can take any positive value.
To compute the update step for $z^{k+1}$ we need to solve a positive-definite
system of equations by only accessing $\mathcal A^*$ and $\mathcal A$.
There are many methods to do this, for example, conjugate-gradient (CG)
type methods \cite{hestenes1952methods,krylov1931numerical,shewchuk1994introduction}.

For all of these algorithms, $b^k$ converges to a solution of \eqref{e-estimators}.
There are many other algorithms for evaluating these estimators; 
see, \eg, \cite{beck2017first,chambolle2011first,nesterov2013gradient,nw2006numopt,simon2013blockwise}.
The methods for computing SURE we describe below will work with most of these as well.

\subsection{Weak differentiability of convex regularized
  regression}\label{s-weak-diff}

For SURE to be an unbiased estimate of risk, the estimator $\hat\mu$
must be weakly differentiable \cite[\S6]{evans15measure} and
obey some integrability conditions \cite{stein1981estimation}.
For our purpose, it is sufficient to show that $\hat\mu$ is Lipschitz continuous
\cite[Lemma III.2]{candes2013svt}.

We will now show that $\hat\mu$ is Lipschitz if $r$ is a closed convex 
proper function.
The coefficient estimate $\hat \beta(y)$ minimizes
$r(b) + \frac{1}{2} \|\mathcal A b - y\|_2^2$,
and so by \cite[Thereom 3.1.23]{nesterov2018lectures}, we have for all $w$,
\begin{multline*}
    r(w) \geq r(\hat\beta(y)) + \langle\mathcal A^* y
        - \mathcal A^* \mathcal A\hat\beta(y)\mid w - \hat\beta(y)\rangle =\\
              r(\hat\beta(y)) + \langle y
        - \mathcal A\hat\beta(y)\mid \mathcal Aw - \mathcal A\hat\beta(y)\rangle
        .
\end{multline*}

Evaluating this at $w =\hat \beta(\tilde y)$ gives
\[
r(\hat\beta(\tilde y)) \geq r(\hat\beta(y)) + \langle y
- \mathcal A\hat\beta(y) \mid \mathcal A\hat\beta(\tilde y) - \mathcal A\hat\beta(y)\rangle,
\]
and switching the roles of $y$ and $\tilde y$, we obtain
\[
r(\hat\beta(y)) \geq r(\hat\beta(\tilde y)) + \langle  \tilde y 
- \mathcal A\hat\beta(\tilde y)
\mid \mathcal A\hat\beta(y) - \mathcal A\hat\beta(\tilde y)\rangle.
\]
Adding these two inequalities yields 
\BEAS
    0 &\geq&
    \langle  y -  \mathcal A\hat\beta(y)\mid \mathcal A  \hat\beta(\tilde y) - \mathcal A \hat\beta(y)\rangle \\ & &\hspace{2em}+ 
\langle \tilde y - \mathcal A\hat\beta(\tilde y) \mid \mathcal A \hat\beta(y) - \mathcal A \hat\beta(\tilde y)\rangle \\
&=& 
    \langle y - \hat\mu(y) - \tilde y + \hat\mu(\tilde y) \mid \hat\mu(\tilde y) - \hat\mu(y)\rangle \\
&=&
    \langle y - \tilde y\mid \hat\mu(\tilde y) - \hat\mu(y)\rangle +  \|\hat\mu(\tilde y) - \hat\mu(y) \|_2^2.
\EEAS
Re-arranging and using the Cauchy-Schwartz inequality gives 
\begin{multline*}
\|\hat\mu(\tilde y) - \hat\mu(y) \|_2^2 \leq
\langle \tilde y - y \mid \hat\mu(\tilde y) - \hat\mu(y)\rangle \\
    \leq \|\tilde y - y\|_2 \|\hat\mu(\tilde y) - \hat\mu(y) \|_2,
\end{multline*}
eliminating a factor of $\|\hat\mu(\tilde y) - \hat\mu(y) \|_2$ shows that
$\hat\mu$ is $1$-Lipschitz.

\section{SURE-CR}\label{s-sure-cr}

\subsection{Randomized trace estimation}\label{s-trace-est}

In this section we describe methods for estimating the trace of a
$d \times d$ matrix $M$, that access $M$ only via an oracle that evaluates its 
adjoint, $v\mapsto M^*v$.
We refer to this oracle as vector-matrix oracle, since it evaluates (the transpose of) 
$v^*M$.
We will apply this to the specific matrix $M=D\hat\mu(y)$ to evaluate 
the divergence term in SURE.

The na\"ive approach is to use the oracle to evaluate $M^*e_i$,
where $e_i$ is the $i$th unit vector, for $i=1, \ldots, d$, whereupon
we can readily evaluate
\[
\Tr M = \sum_{i=1}^d e_i^* (M^* e_i).
\]
When $d$ is very large, this is slow. It also evidently involves much wasted 
computation, since we end up computing all $d^2$ entries of $M$, 
only to sum the $d$ diagonal ones.

Randomized methods can be used to estimate $\Tr M$ using far fewer than $d$ 
evaluations of the adjoint mapping.
These methods are based on the simple observation that if the random 
variable $Z\in \reals^d$ satisfies $\Expect Z =0$ and $\Expect ZZ^* = I$,
then we have $\Expect Z^*MZ = \Tr M$.
To approximate this we compute $m$ independent samples of $Z$, 
$z_1, \ldots, z_m$, and take the empirical mean as our estimate,
\[
\Tr M \approx \frac{1}{m}\sum_{i=1}^m z_i^* (M^* z_i),
\]
which is unbiased.
In \cite{hutchinson1989}, Hutchinson showed that the variance of the error 
in this approximation is minimized if the $Z_i$'s are i.i.d.~random variables
taking values $\pm 1$, each with probability $1/2$, which is known as the
Rademacher distribution.

Improvements on this basic randomized method
were recently suggested by Meyer, Musco, Musco, and Woodruff
in \cite{meyer2021hutch++}.
They proposed \hutchpp, which uses a low-rank approximation of the matrix to
project some queries away from large singular values of the matrix.
\hutchpp\ is also an unbiased estimator of the trace, and consistently produces 
a good estimate of the trace using fewer
queries to the vector-matrix oracle than the basic randomized method.
\hutchpp's computation takes part in three phases, each of which requires an equal
number of calls to the vector-matrix oracle, so the total number of queries is a
multiple of $3$.
In the first phase, \hutchpp\ sketches the matrix;
\ie, it multiplies $M$ with a 
tall rectangular matrix whose entries are i.i.d.~Rademacher random variables,
and computes an orthogonalization of that matrix product, which is an estimate of
the dominant dimensions of the matrix.
In the second phase, it computes the exact trace of $M$ projected onto the
dominant dimensions found via the sketch.
In the third phase, it runs the Hutchinson estimator on $M$ projected away
from those dominant dimensions.

In our method for evaluating SURE, we found that $34$ queries per
\hutchpp\ phase, for a total of $102$ vector-matrix oracle calls, consistently
produced high quality estimates of the trace.
For small problems, \ie, those of size less than or equal to $102$, we
exactly compute the trace without any randomization.

Subsequent works have developed alternative trace estimation algorithms
\cite{persson2022,epperly2023xtrace}.

\subsection{Vector-Jacobian oracles}\label{s-vec-jac}

In this section we describe methods for computing the adjoint oracle 
$v \mapsto \left(D\hat \mu(y)\right)^* v$.
Using $\hat \mu(y) = \mathcal A \hat\beta(y)$, we have 
$D\hat \mu(y)=\mathcal A D\hat\beta (y)$ and, therefore, 
\[
    \left(D\hat\mu(y)\right)^* v =
\left(D \hat\beta(y)\right)^* (\mathcal A^* v).
\]
So it suffices to evaluate the mapping $u \mapsto\left( D\hat\beta(y)\right)^* u$.
Roughly speaking, we need to differentiate through the solution of 
the optimization problem \eqref{e-estimators}, \ie, the mapping from the data
$y$ to the parameter estimate $\hat\beta(y)$.

\paragraph{Differentiability}
In many cases $\hat\mu$ is not differentiable.
However in \S\ref{s-weak-diff}
we showed that $\hat\mu$ is Lipschitz;
by applying Rademacher's theorem, 
we know that $\hat\mu$ is a.e.-differentiable
under the Lebesgue measure, and since $y$ has a Gaussian 
distribution $\hat\mu$
is almost surely differentiable at $y$ \cite[\S3.1.2]{evans15measure}.

\paragraph{Generic methods}

Some recent work shows how to differentiate through the solution of some 
convex optimization problems
(when the mapping is differentiable), for example \cite{amos2017optnet} for quadratic programs 
(QPs) and \cite{diffcp2019} for cone programs.
These methods in turn have been integrated into
software frameworks for automatic differentiation such as PyTorch \cite{pytorch} 
and TensorFlow \cite{tensorflow2015-whitepaper,agrawal2019tensorflow}.
Such libraries include CVXPYlayers, diffcp, and OPTNET
\cite{cvxpylayers2019,diffcp2019,amos2017optnet}.
All of these give methods for evaluating $u \mapsto D\hat\beta(y)^* u$, without 
forming the matrix $D\hat\beta(y)$.
These generic methods work well for small problems and some medium-sized
problems, but they do not scale to large scale problems.
At non-differentiable points, these methods compute a heuristic quantity
\cite[\S14]{griewank2008evaluating}.

\paragraph{Differentiating through an iterative solver}

Another approach to differentiating through a convex problem relies on a solver
or iterative solution algorithm, such as those described in \S\ref{s-solver-algorithms}.
Existing work differentiates through proximal operators to use them as non-linear
activations in neural networks \cite{diamond2018unrolled,wang2016proximal}, in this
work, we differentiate through iterative optimization algorithms to approximate
differentiating the solution map.
Here we view the iterative algorithm as a sequence of mappings, \ie,
we view our iterative algorithm as applying an operator $F^k$ at each iteration such that
\[
    b^{k+1}, S^{k+1} = F^k(b^k, S^k, y)
\]
where $S^k$ is any ancillary state in the algorithm
(\eg, in FISTA $S^k = b^{k-1}$ and in ADMM $S^k = (z^{k}, u^k)$).
Suppose it takes $\ell$ iterations to converge to a reasonable tolerance,
so $\hat\mu(y) \approx F^\ell(F^{\ell -1}(\ldots, y), y)$.
By implicitly differentiating this recurrence and applying the chain rule,
we obtain a series of equations that we use to compute
$(D \hat\mu(y))^* v$, given a vector of output sensitivities $v$.
Our approximation of $\hat\mu$ may be non-differentiable on a set of positive
Lebesgue measure.
In this situation, we need to compute a quantity that can serve
as a surrogate for the true vector-Jacobian product.
In neural network training, it is common to discuss the vector-Jacobian of a scalar
loss function---which is simply the gradient---even when the loss function is
non-differentiable.
Many choices of surrogates for when differentiability fails 
have been proposed and seem to work well here \cite{begnio2013estimating,glorot2011deep}.
In \S\ref{s-matrix-compl-example}, our empirical results show that a continuous
extension of the true derivative yields sufficiently accurate estimates at
non-differentiable points of the derivative of $\hat\mu$ so that we still have a
good estimate of the risk of $\hat\mu$.

As an example of differentiating our approximation of $\hat\mu$, we work through
the derivative of ISTA.
ISTA is straightforward to analyze because there is no ancillary state in the
algorithm, but this method easily generalizes to the other algorithms from
\S\ref{s-solver-algorithms}.
To simplify our equations, we let
$b^{k+1/2} = b^{k} -\eta \mathcal A^*\left( \mathcal A b^{k} - y\right)$.
By differentiating the ISTA iterations we obtain
\[
    \begin{array}{ll}
        Db^{k+1} &= D\prox_{\eta r}\left(b^{k+1/2}\right) Db^{k}
                        - \eta \mathcal A^* \left( \mathcal A Db^{k}
                        -I\right) \\
                    &= \left(D\prox_{\eta r}\left(b^{k+1/2}\right)
                        -\eta \mathcal A^* \mathcal A\right) Db^{k}
                        + \eta \mathcal A^* 
                        .
    \end{array}
\]

In a forward pass, we can evaluate $b^{k+1/2}$  for $k = 1, \ldots, \ell - 1$
and cache them to enable the vector-Jacobian oracle evaluations.
Evaluating $(Db^{\ell})^* v$ then becomes a recursive problem, which can be
computed using two of the oracles we needed for the forward pass---$\mathcal A$
and $\mathcal A^*$---and one new oracle: the vector-Jacobian oracle for the proximal
operator.
The base case for our recursion comes from our requirement that $b^1$ is chosen
independently of $y$ \ie\ that $(Db^1)^* v = 0$.
The difficulty in this method relies in evaluating the vector-Jacobian oracle of
the proximal operator.

\paragraph{Evaluating proximal operator vector-Jacobian oracles}

For many proximal operators known in closed-form, the Jacobians are trivial to find
in closed-form.
For example, the $\ell_1$ norm has proximal operator given by soft-thresholding,
$\mathcal T_\eta$.
Since soft-thresholding occurs component-wise, this means that the Jacobian is a
diagonal matrix, whose non-zero entries are
$1$ if $b^{k+1/2}_i$ is above the threshold, $-1$ if it is below the negative
of the threshold, and $0$ otherwise.
In this case, it is possible to efficiently compute the vector-Jacobian oracle
without forming the whole Jacobian to find
\[
    \left(D \prox_{\eta \|\cdot\|_1}\left(b^{k+1/2}\right)\right)^* u =
    \diag\left(J_{\mathcal T_\eta}\left(b^{k+1/2}\right)\right) \circ u ,
\]
where $a \circ b$ denotes Hadamard or component-wise multiplication.
Here we handle the points of non-differentiability by using the value of the
derivative at a point in a very small neighborhood of the non-differentiable
point.
In particular, since non-differentiability occurs only when entries of $X^Ty$ are
exactly
equal to $\eta$, we can instead interpret our choice of a point in the neighborhood
as evaluating the derivative at a point within the floating point uncertainty of
our vector.

Other closed-form proximal operators have non-trivial Jacobians.
For example, the proximal operator of the nuclear norm is given by
\[
    \prox_{\eta \|\cdot\|_*}(b^{k+1/2}) = U \mathcal T_\eta(\Sigma) V^*,
\]
where $b^{k+1/2} = U \Sigma V^*$ is the singular-value decomposition of
$b^{k+1/2}$ and $\mathcal T_\eta$ is soft-thresholding on $\Sigma$.
This has a non-trivial Jacobian because of the multi-valued nature of the SVD
in the presence of repeated singular values.
However, since all proximal operators are Lipschitz, we know that it is
a.e.-differentiable.
\cite[Lemma IV.2]{candes2013svt} gives closed-form expressions for the
Jacobian of this proximal operator that hold for simple and full-rank matrices.
However, it is common that later iterations will involve low-rank matrices, which
requires us to select an approximation of the vector-Jacobian products.
We use the continuous extension of the closed-form vector-Jacobian product,
which exists for all matrices which do not have any singular values exactly equal
to $\eta$.
We derive an expression for this extension in \S\ref{s-grad-prox-nuc-norm}.
For matrices with singular values exactly equal to $\eta$ we just evaluate at a
point within the neighborhood of the matrix similar to how we handle the $\ell_1$
norm.

In general, when trying to apply SURE-CR to a new proximal operator, it is necessary
to be able to evaluate the vector-Jacobian product for that proximal operator.
If the proximal operator has points of non-differentiability which are reached by
the iterative algorithm, then it is necessary to choose a surrogate for the
vector-Jacobian product.
The accuracy of SURE-CR is limited by the accuracy of the 
vector-Jacobian product oracle.

Since $\left(D\hat\mu(y)\right)^* v = \nabla_y \langle\hat\mu(y) \mid v\rangle$,
it is possible to apply well-known strategies to compute the gradient of a
scalar-valued function.
Most notably, reverse-mode automatic differentiation automates much of this
section's work \cite{griewank1989automatic}.
For many proximal operators with closed-form expressions, reverse-mode automatic
differentiation can differentiate the proximal operator without an analytic
derivation of a closed-form for the vector-Jacobian oracle.

As an example we work out how to construct the oracle for $r(b) = \|b\|_1 + \|b\|_2^2$
(a weighted sum of these norms is the regularizer in the elastic
net \cite{zou2005regularization}).
The proximal operator can be evaluated by applying separability to find that 
it is given by a scaled form of soft-thresholding,
\BEAS
    \prox_{\eta r}(v) &=& \argmin_b \left(  \eta \|b\|_1 + \eta \|b\|_2^2 +
    \frac{1}{2}\|b - v\|_2^2 \right) \\ &=& \frac{1}{1 + 2 \eta}
 \mathcal T_\eta(v) .
\EEAS
By rewriting soft-thresholding as 
\[
\mathcal T_\eta(v) = (v - \eta\ones)_+ - (-v - \eta\ones)_+,
\]
we can express this function in terms of elementary operations that are commonly
supported by automatic differentiation libraries, meaning no work is required to
construct the vector-Jacobian oracle.

\subsection{Implementation}\label{s-implementation}

We have implemented the methods described above in SURE-CR,
an open-source package available at
\begin{center}
\url{https://github.com/cvxgrp/SURE-CR}.
\end{center}
It supports divergence computation via CVXPYlayers as well as 
via differentiation through FISTA and ADMM, and uses \hutchpp\ to 
estimate the divergence.

SURE-CR relies on an existing computational graph library, pyTorch
\cite{pytorch}, to enable GPU-acceleration in our solvers and to enable
reverse-mode automatic differentiation.
We have implemented a library to encode the linear operator $\mathcal A$ as
a computational flow graph.
It is available at \begin{center}
\url{https://github.com/cvxgrp/torch\_linops}.
\end{center}
This library adapts Barratt's preconditioned conjugate gradient 
implementation \cite{barratt2019}
and implements randomized preconditioners including Nystr\"om preconditioning
\cite{frangella2021randomized}.

By differentiating through FISTA and ADMM iterations, SURE-CR is able to scale
to large problems.
For example it can evaluate SURE for a matrix completion problem with 
$b \in \reals^{2000 \times 1000}$ and $10\%$ of entries revealed, for which
$D\hat \mu(y)$ is a $10^5 \times 10^5$ matrix (which of course is never formed)
in $120$ seconds on the server described in \S\ref{s-examples}.

To apply SURE-CR to novel problems and regularizers, the user should adapt an
example from \S\ref{s-code-examples} by implementing their linear operator
$\mathcal A$ and $\mathcal A^*$ as shown in \S\ref{s-mat-compl-code-example}
and implementing the proximal operator as a differentiable torch function.
This can be done most easily by expressing it as the composition of built-in
torch functions as shown in \S\ref{s-lasso-code-example}.
In the event that a heuristic is used for the derivative of the proximal operator,
it may be valuable to test that the heuristic and the true vector-Jacobian
products found by CVXPYlayers agree.

SURE-CR currently uses at most one GPU; however, in
hyperparameter sweep problems, users can run different experiments on different 
GPUs in parallel.

\section{Numerical examples}\label{s-examples}

In this section we report results of numerical examples of SURE-CR.
We consider three problems, LASSO, matrix completion, and robust PCA, 
and for each one, problem instances ranging from small to large.
For each instance we evaluate various estimates of SURE, as well as 
an estimate of the $\ell_2$ risk obtained via a Monte Carlo method described
below.

We carry out a few additional experiments that analyze the 
variance contributed by
SURE itself in high-dimensions, and also, the variance contributed by our
use of a randomized trace estimator.
We will see that the latter is substantially smaller than the former.

Finally, in our last example, we show how SURE-CR can be used to carry 
out hyperparameter selection.

\paragraph{Hyperparameter selection}
When selecting regularization parameters, we swept over the parameters---equally
spaced on a logarithmic scale---on the largest problem size we planned to run.
We then selected a value which had
risk less than half the risk of the maximum likelihood estimator of $\mu$ and
had a high iteration count relative to the other runs in the sweep.
We require the risk to be small in order to demonstrate SURE-CR in problem settings
where the estimator is useful.
The higher the iteration count, the longer SURE-CR takes to run since we have
to differentiate through more iterations of the solver algorithm; accordingly, to
give a better sense of worst-case runtime when using SURE-CR we prefer problem
instances that gave higher iteration counts.

\paragraph{SURE estimates}
In our first example, LASSO, we report the value of the analytical expression for SURE.
In all examples we evaluate SURE using CVXPYlayers, where it was possible,
\ie, for the smaller problem instances.
For each problem we use either ADMM or FISTA, depending on which 
was faster on small test problems.

\paragraph{Monte Carlo $\ell_2$ risk estimate}
Since we are using
synthetic data and know $\mu = \mathcal A \beta$, we are able to use a Monte Carlo
method to approximate the risk as
\[
    R(\hat\mu) \approx \frac{1}{m}\sum_{i=1}^m \|\hat\mu(y_i) - \mu\|_2^2,
\]
where $y_i \overset{\text{i.i.d.}}{\sim} \mathcal N(\mu, \sigma^2 I)$.
(In practical problem settings, this Monte Carlo estimation of $\ell_2$ risk is
not possible.)

\paragraph{Computational platform}
We report timings for running SURE-CR on the Stanford University Institute
for Computational and Mathematical Engineering's DGX-1, with 8
Nvidia Tesla V100-SXM2-32GB-LS GPUs, an Intel Xeon E5-2698 v4 with 80 cores,
540GiB of memory, and 32GiB of GPU memory per GPU.
(However, we were limited to only one GPU during our tests.)

\paragraph{Overview of results}
The results are summarized in the tables below. 
Comparing the values of the various estimates of SURE and $\ell_2$ risk across
each row, we see that there is good agreement, except for the smallest problem
instances.
In \S\ref{variance-bound}, we show that recent works by Bellec and Zhang
\cite{bellec2021debiasing,bellec2021second} enables bounding the variance of SURE
to be less than $4\sigma^4 d + 2 \sigma^2 R(\hat\mu)$.
For our estimators, the risk scales about affinely with $d$, and therefore
the standard deviation of SURE grows slower than its expectation, so we see
asymptotic convergence to the true value in relative error.

For the largest instances of matrix completion and 
robust PCA, each of which have 2 million parameters, we are able to compute
SURE in under two minutes.
To our knowledge, there was no previously known method for computing SURE for
such large instances.

\subsection{LASSO}\label{s-lasso-example}
We compute SURE for LASSO problems, described in \S\ref{s-classical-estimators}.
We consider under-determined problems with $p = 2d$, for
$d = 250, 500, 2500, 5000, 25000$.

\paragraph{Data generation}
We draw the entries of the data matrix i.i.d.~from a standard normal
distribution.  We pick $\beta$ with $d/20$ nonzero entries equal to
a constant and use $\sigma^2 = 2$.
We pick the value of the nonzero coefficients so that
$\frac{\|\mu\|_2^2}{\|\mu\|_2^2 + d\sigma^2} = 0.8$.
We sample one $y$ independently from the rest of the data and select
$\lambda = 0.1\lambda_\text{max}$ (defined in \eqref{e-lambda-max}).

For each instance we use SURE-CR with CVXPYlayers,
SURE-CR with FISTA, the analytic SURE value computed against CVXPY's solution,
and the Monte Carlo estimate of the risk using CVXPY to solve
the optimization problem.
In its default configuration, CVXPYlayers has very low accuracy in moderate
dimensions and does not raise warnings about the errors. 
To correct for this, we switched CVXPYlayer's implicit linear system solver for
its direct linear system solver; this did not significantly impact runtime on
problems for which it was giving accurate results.
We present both the risk and time values for each.
When using CVXPYlayers, we report the value as $*$ when CVXPYlayers has a
non-standard return status warning, raises an error, or takes more than 12 hours.
The seed used to generate the \hutchpp\ queries and the sample point at which to
compute SURE are the same for all problems of a given size.
The results are given in Table~\ref{t-lasso}.

\begin{table*}
\begin{center}
\caption{Values of coordinate-wise SURE estimates and computation times for
    five LASSO problem instances.
    Coordinate-wise SURE is given by $(1/d){\SURE(\hat\mu, y)}$ and is used to
    improve readability.
    Times given in seconds.}
\label{t-lasso}
\begin{tabular}{ll|ll|ll|ll}
    \toprule 
\multicolumn{2}{c}{Dimensions} &
\multicolumn{2}{c}{CVXPYlayers} &
\multicolumn{2}{c}{FISTA} &
Analytic & MC risk\\
$d$ & $p$ & Value & Time & Value & Time \\
    \midrule 
       250 &    500 & 0.51 & 252    & 0.51 & 2.54 & 0.52 & 0.48(\textless 0.005) \\
       500 &   1000 & 0.43 & 1929   & 0.37 & 2.89 & 0.40 & 0.50(\textless 0.005) \\
      2500 &   5000 & *    & *      & 0.58 & 3.10 & 0.57 & 0.51(\textless 0.005) \\
      5000 &  10000 & *    & *      & 0.47 & 9.19 & 0.46 & 0.53(\textless 0.005) \\
     25000 &  50000 & *    & *      & 0.58 & 287  & 0.58 & 0.54(\textless 0.005) \\
    \bottomrule 
\end{tabular}
\end{center}
\end{table*}

\subsection{Matrix completion}\label{s-matrix-compl-example}
\begin{figure}
    \resizebox{\columnwidth}{!}{\input{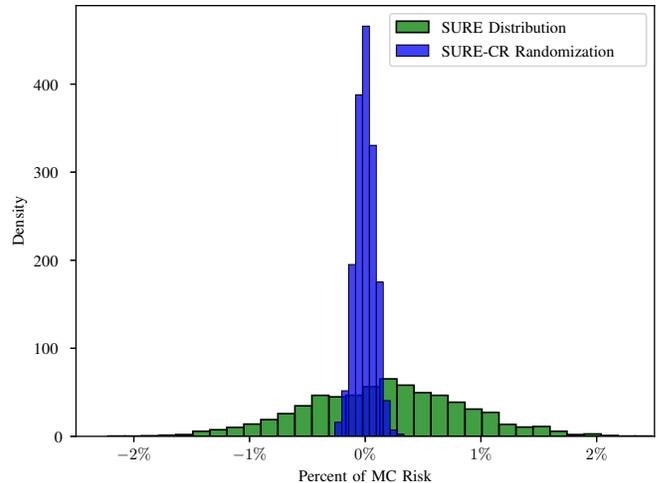}}
    \caption{The green histogram is the relative error between SURE at various 
    sample points against the Monte Carlo risk.
    The blue histogram shows the relative error between SURE-CR at a sample point 
    and the mean of $100$ runs of SURE-CR at that point.\label{f-matrix-compl-histogram}}
\end{figure}

We compute SURE for matrix completion problems, described in
\S\ref{s-matrix-estimators}.  For all problems, we use $d = 0.1 m n$,
$\sigma^2 = 2$, and $\lambda = 0.25 \lambda_{\max}$.  For the large
problems used to generate Figure~\ref{f-matrix-compl-histogram}, we
use $m = 2000$, $n = 1000$, $d = 0.1 m n = 200000$. 
  We use SURE-CR with CVXPYlayers and SURE-CR with
ADMM to compute SURE in Table~\ref{t-matrix-compl}.
We describe how we formed $\mu$ and $\beta$ below.

\begin{table}
\begin{center}
\caption{Values of coordinate-wise SURE estimates and computation times for five matrix
    completion problem instances.
    Coordinate-wise SURE is given by $(1/d){\SURE(\hat\mu, y)}$ and is used to
    improve readability.
    Times given in seconds.\label{t-matrix-compl}}
\begin{tabular}{ll|ll|ll|l}
    \toprule 
\multicolumn{2}{c}{Dimensions} &
\multicolumn{2}{c}{CVXPYlayers} &
\multicolumn{2}{c}{ADMM} &
 MC risk\\
$d$ & $p$ & Value & Time & Value & Time \\
    \midrule 
    20             & 200            & 1.15 & 1.51  & 1.16 & 5.20 & 1.33(0.01) \\
    500            & 5000           & 0.86 & 2246  & 0.88 & 49.9 & 0.96(\textless 0.005) \\
    2000           & $2\times 10^4$ & 0.84 & 15866 & 0.84 & 45.1 & 0.90(\textless0.005) \\
    $5\times 10^4$ & $5\times 10^5$ & *    & *     & 1.69 & 41.2 & 1.70(\textless0.005) \\
    $2\times 10^5$ & $2\times 10^6$ & *    & *     & 0.74  & 114  & 0.74(\textless0.005) \\
    \bottomrule 
\end{tabular}
\end{center}
\end{table}

Since $\mathcal A^* \mathcal A + \lambda I$ is a diagonal matrix, we replace
the preconditioned conjugate gradient step of the ADMM updates with an exact
inverse.
This has no significant impact on the numerical accuracy of our algorithm, but
does improve its runtime.

\paragraph{Data generation}
We first generate $\beta = U\Sigma V^* \in \reals^{m \times n}$ with 
$\max(5, 0.02 n)$ non-zero singular values, which are uniformly distributed over $[0, n]$.
The matrices $U$ and $V^*$ are generated by computing the SVD of a matrix where each entry
is independent and identically distributed as uniform over $[0, 1]$.
For the selection operator $\mathcal A$, we selected $10\%$ of the entries at
random without replacement.
We then sampled $y \sim \mathcal N(\mathcal A \beta, \sigma^2 I)$.

\paragraph{Quantifying \hutchpp\ uncertainty}
We verify that the uncertainty from using \hutchpp\ to estimate the
divergence is dominated by the uncertainty inherent in SURE.
For $20$ sample points of $y$, we ran SURE-CR on each point $100$ times.
In Figure~\ref{f-matrix-compl-histogram}, we show in blue the distribution of
the relative error between the SURE-CR values and the sample mean of the SURE-CR
runs on that point: let  $\text{SURE-CR}(y, i)$ denote the
random variable of the output of running SURE-CR on a point 
$y$ with seed $i$.
Then for samples $y_1, y_2, \ldots, y_{10\,020}$, we plot the histogram of 
\[
    \frac{
        \begin{array}{l}
            \text{SURE-CR}(\hat\mu, y_i, 100 i + j) \hspace{10em}\mbox{} \\ \hfill -
        100^{-1}\sum_{k=1}^{100}\text{SURE-CR}(\hat\mu, y_i, 100 i + k)
        \end{array}
    }{
        10\,000^{-1}\sum_{k=21}^{10\,020}
        \|\hat\mu(y_k) - \mu \|_2^2
    }
\]
for $i=1,2, \ldots, 20$ and $j = 1, 2, \ldots, 100$.
We also plot the histogram of the relative error between $2000$ evaluations of
SURE-CR and the Monte Carlo estimation of the risk.
The uncertainty from the algorithm's randomization is small compared to SURE's
uncertainty.

\paragraph{SURE as estimate of risk}
The green histogram in Figure~\ref{f-matrix-compl-histogram} shows that SURE-CR is
within $2.5\%$ of the Monte Carlo risk at $2000$ independent sample points.
Precisely, the green histogram shows the histogram of the quantity
\[
    \frac{
        \text{SURE-CR}(\hat\mu, y_i, i) - 10\,000^{-1}\sum_{j=2001}^{12\,000}\| \hat\mu(y_j) - \mu \|_2^2
    }{
        10\,000^{-1}\sum_{j=2001}^{12\,000}\| \hat\mu(y_j) - \mu \|_2^2
    }
\]
for $i =1, 2, \ldots, 2000$ and independent samples $y_1, y_2, \ldots, y_{12\,000}$.
This shows SURE is a good estimate of the true risk.

\paragraph{Non-differentiability}
In around $5\%$ of the 2000 samples used to generate the green histogram in
Figure~\ref{f-matrix-compl-histogram}, we observed that our approximation of
$\hat\mu$ was non-differentiable.
We detected this by running our algorithm without using the extension of the derivative
and seeing what percentage of runs encountered numerical issues caused by repeated 
or zero singular values.
We then ran the experiment using the extension of the derivative, and report those
values here.
Notably, those samples are indistinguishable from the other samples in the
histogram, showing that our heuristic is effective at approximating the
vector-Jacobian products for $\hat\mu$ and still providing a good estimate of risk.

\subsection{Robust PCA}\label{s-robust-pca-example}
\begin{figure}
    \resizebox{\columnwidth}{!}{\input{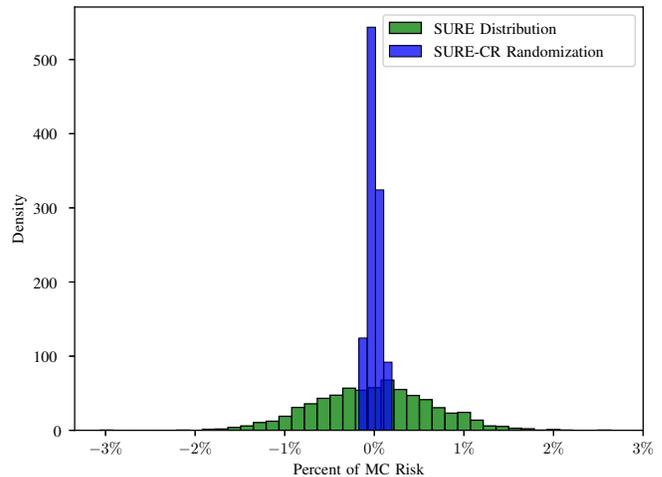}}
    \caption{The green histogram is the relative error between SURE at various 
    sample points against the Monte Carlo risk.
    The blue histogram shows the relative error between SURE-CR at a sample point 
    and the mean of $100$ runs of SURE-CR at that point.\label{f-robust-pca-histogram}}
\end{figure}

We also tested SURE on robust PCA problems, described in
\S\ref{s-matrix-estimators}.
For all problems, we use $m = n$, $\sigma^2 = 2$, $\lambda = 0.16\lambda_\text{max}$,
and $\gamma = 0.057\gamma_\text{max}$.
For the large problems used to generate Figure~\ref{f-robust-pca-histogram}, we
used $m = n = 1000$.
We use SURE-CR with CVXPYlayers and SURE-CR with ADMM to compute SURE in
Table~\ref{t-robust-pca}.

\begin{table}
\begin{center}
    \caption{Values of coordinate-wise SURE estimates 
    and computation times for five robust PCA problem instances.
    Coordinate-wise SURE is given by $(1/d){\SURE(\hat\mu, y)}$ and is used to
    improve readability.
    Times given in seconds.\label{t-robust-pca}}
\begin{tabular}{ll|ll|ll|l}
    \toprule 
\multicolumn{2}{c}{Dimensions} &
\multicolumn{2}{c}{CVXPYlayers} &
\multicolumn{2}{c}{ADMM} &
 MC risk\\
$d$ & $p$ & Value & Time & Value & Time \\
    \midrule 
    100              & 200            & 5.14 & 1.15 & 5.13 & 16.0 & 5.01(0.006) \\ 
    2500             & 5000           & 0.53 & 115  & 0.53 & 19.9 & 0.59(\textless0.005) \\ 
    10000            & $2\times 10^4$ & 0.31 & 1116 & 0.31 & 21.5 & 0.34(\textless0.005) \\
    $2.5\times 10^5$ & $5\times 10^5$ & *    & *    & 0.27 & 22.1 & 0.27(\textless0.005) \\
    $1\times 10^6$   & $2\times 10^6$ & *    & *    & 0.44  & 31.4 & 0.44(\textless0.005) \\
    \bottomrule 
\end{tabular}
\end{center}
\end{table}

\paragraph{Data generation}
We select $S$ with $\max(10, 10^{-4} n^2)$ non-zero entries drawn from a uniform
distribution over $[0, 100]$.
We select $L$ with rank $\max(5, 0.02 n)$ and singular values distributed uniformly
over $[0, n]$.
We sampled $y \sim \mathcal N(L + S, \sigma^2 I)$.

\paragraph{SURE as estimate of risk}
Figure~\ref{f-robust-pca-histogram} shows the histogram of the relative error
compared to the Monte Carlo estimate of the risk for $m = n = 1000$ and
the histogram of the variance from the randomization in SURE-CR for $m = n = 1000$.
We ran SURE-CR on $2000$ sample points and use $10\,000$ samples for the Monte Carlo estimate.
We observed only one sample for which SURE-CR diverged from the Monte Carlo risk by more than $3\%$.

\subsection{SURE for hyperparameter selection}
\begin{figure}
    \resizebox{\columnwidth}{!}{\input{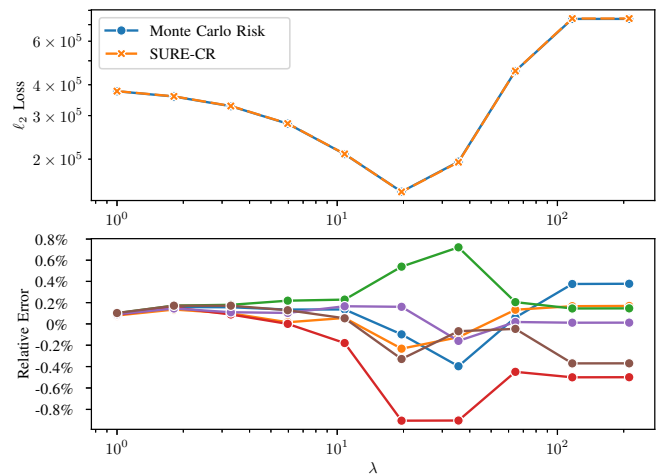}}
    \caption{\emph{Top}. SURE-CR and Monte Carlo estimate of $\ell_2$ risk
as a function of the hyperparameter.
A single sample of $y$ was used for all of the SURE-CR runs.
The two lines are visually indistinguishable.
    \emph{Bottom}.  Relative error plots for the SURE-CR sweep run on 6 independent
    samples of $y$.  The Monte Carlo estimate and the computed SURE value
differ by less than $1\%$.\label{f-hyperparameter-search}}
\end{figure}
In this experiment, we aim to select an optimal hyperparameter for matrix
completion.
We use the same setup as in \S\ref{s-matrix-compl-example} with $m = 2000$ and $n = 1000$,
except we now draw a single sample $y$.

We then ran a grid search over $\lambda$, varying it exponentially over
$[1, 2 \lambda_{\max}]$, where $\lambda_{\max}$ is the smallest $\lambda$ for
which $\hat\beta(y) = 0$.
We drew a single sample of $y$, and then for each $\lambda$ we ran SURE-CR with ADMM.
We then computed a Monte Carlo estimation of the risk for each $\lambda$.
Figure~\ref{f-hyperparameter-search}, shows that the risk versus $\lambda$ curves
are visually indistinguishable.
We also show that for 6 independent samples of $y$, the relative error was
consistently below $0.9\%$.

\section*{Acknowledgments}

We thank Mert Pilanci for many helpful comments during a talk about this
project.  We thank Raphael Meyer for help with \hutchpp. We also thank an
anonymous reviewer for an unusually thorough and careful review that helped
improve the paper.

\bibliography{sure}

\appendices

\section{Code examples}\label{s-code-examples}
\subsection{CVXPYlayers --- LASSO}
This code sample demonstrates how to use SURE-CR with CVXPYlayers.
It is based on the code used in \S\ref{s-lasso-example}.
\begin{lstlisting}[language=mypython]
import cvxpy as cp

import surecr

X, y, variance, lambda_val = ...  # Generate data

beta_cvx = cp.Variable(X.shape[1])
y_cvx = cp.Parameter(y.shape[0])
prob = cp.Problem(
        cp.Minimize(
            1/ 2 * cp.sum_squares(X @ beta_cvx - y_cvx)
            + lambda_val * cp.pnorm(beta_cvx, 1)
        ))
solver = surecr.CVXPYSolver(prob, y_cvx, [beta_cvx], lambda b: X @ b)
sure = surecr.SURE(variance, solver)
cvx_sure_val = sure.compute(y)
\end{lstlisting}

\subsection{LASSO}\label{s-lasso-code-example}
This code sample demonstrates how to use SURE-CR with FISTA and how to define a
custom proximal operator.
It is based on the code used in \S\ref{s-lasso-example}.
\begin{lstlisting}[language=mypython]
import torch
import surecr
import linops as lo

X, y, variance, lambda_val = ... # Generate data

d, p = X.shape
A = lo.aslinearoperator(X.cuda())
y_cuda = y.cuda()
def prox(v, t):
    return torch.relu(v - lambda_val * t) - torch.relu(-v - lambda_val * t)

solver = surecr.FISTASolver(
        A, prox, torch.zeros(p).cuda(),
        device=y_cuda.device)
sure = surecr.SURE(variance, solver)

sure_val = sure.compute(y_cuda)
\end{lstlisting}

\subsection{Matrix completion}\label{s-mat-compl-code-example}
This code sample demonstrates how to use SURE-CR with ADMM and how to define a
custom linear operator.
It is based on the code used in \S\ref{s-matrix-compl-example}.
\begin{lstlisting}[language=mypython]
import torch
import surecr
import surecr.prox_lib as pl
import linops as lo

revealed_indices, y, variance, lambda_val, m, n = ... # Generate data


class SelectionOperator(lo.LinearOperator):
    def __init__(self, shape, idxs):
        self._shape = shape
        self._adjoint = _AdjointSelectionOperator(idxs,
                (self._shape[1], self._shape[0]), self)
        self._idxs = idxs

    def _matmul_impl(self, X):
        return X[self._idxs]

    def solve_I_p_lambda_AT_A_x_eq_b(self, lambda_, b):
        LHS = torch.ones_like(b)
        LHS[self._idxs] += lambda_
        return b / LHS

class _AdjointSelectionOperator(lo.LinearOperator):
    def __init__(self,  idxs, shape, adjoint):
        self._shape = shape
        self._adjoint = adjoint
        self._idxs = idxs

    def _matmul_impl(self, y):
        z = torch.zeros(self.shape[0], dtype=y.dtype, device=y.device)
        z[self._idxs] = y
        return z.reshape(-1)

d = len(revealed_indices)
p = m * n
A = SelectionOperator((d, p), revealed_indices)
y_cuda = y.cuda()

# Generates a function that applies singular value thresholding, which uses a
# continous extension of the derivative for the .backward method.
prox = pl.make_scaled_prox_nuc_norm((m, n), lambda_val)

solver = surecr.ADMMSolver(A, prox, torch.zeros(p).cuda(), device=y_cuda.device)
sure = surecr.SURE(variance, solver)

sure_val = sure.compute(y_cuda)
\end{lstlisting}

\subsection{Robust PCA}
This code sample demonstrates how to use SURE-CR with ADMM and how to use advanced
features of torch\_linops to generate the linear operator.
It is based on the code used in \S\ref{s-robust-pca-example}.
\begin{lstlisting}[language=mypython]
import torch
import surecr
import surecr.prox_lib as pl
import linops as lo

y, variance, lambda_val, gamma_val, m, n = ... # Generate data

d = m * n
p = 2 * d
A = lo.hstack([lo.IdentityOperator(d), lo.IdentityOperator(d)])
y_cuda = y.cuda()

# Generates a function that applies singular value thresholding, which uses a
# continous extension of the derivative for the .backward method.
prox_L = pl.make_scaled_prox_nuc_norm((m, n), lambda_val)

def prox_S(v, t):
    return torch.relu(v - gamma_val * t) - torch.relu(-v - gamma_val * t)

def prox(v, t):
    return torch.hstack([
            prox_L(v[:d], t), prox_S(v[d:], t)
    ])

solver = surecr.ADMMSolver(A, prox, torch.zeros(p).cuda(), device=y_cuda.device)
sure = surecr.SURE(variance, solver)

sure_val = sure.compute(y_cuda)
\end{lstlisting}

\section{Bound on the variance of SURE}\label{variance-bound}

In \cite[Theorem 3.2]{bellec2021second}, it is shown that for convex
regularized regression
\[
    \var(\SURE(\hat\mu, y)) \leq
        \Expect[(\SURE(\hat\mu, y) - \|\hat\mu(y) - \mu\|_2^2)^2] + \sigma^4 d
\]
and
\[
    (\SURE(\hat\mu, y) - \|\hat\mu(y) - \mu\|_2^2)^2 \leq
        2 \sigma^2(\|y - \hat\mu(y)\|_2^2 + \SURE(\hat\mu, y))
\]
almost surely.
By applying algebraic manipulation and SURE's unbiasedness, we can find that
\[
    \var(\SURE(\hat\mu, y)) \leq
        3 \sigma^4 d - 4\sigma^4 \Expect[\Div \hat\mu(y)] + 4 \sigma^2 R(\hat\mu).
\]

In \cite[Proposition 5.3]{bellec2021debiasing}, it is shown that $D\hat\mu(y)$ is
almost surely positive semi-definite.
This suggests that $\Div \hat\mu(y) = \Tr(D\hat\mu(y)) \geq 0$ almost surely and
lets us conclude that 
\[
    \var(\SURE(\hat\mu, y)) \leq 3 \sigma^4 d + 4 \sigma^2 R(\hat\mu).
\]

\section{Differentiating the proximal operator of the nuclear norm}\label{s-grad-prox-nuc-norm}

The proximal operator of the nuclear norm is given by a spectral function  $F(X)$
such that $F(X) = U F(\Sigma) V^T$ where $U, \Sigma, V^T$ are the full SVD of $X$
and where $F(\Sigma)$ applies the function $\mathcal T_\eta(\sigma) =
(\sigma - \eta)_+$ elementwise to all entries of $\Sigma$.

The function is non-differentiable when $X$ has repeated singular values,
any singular values equal to $0$, or any singular values equal to $\eta$.
Formally, the mapping $X \mapsto \left(DF(X)\right)^* Z$ for a fixed matrix $Z$,
is only defined when $X$ has all distinct singular values and no singular values
equal to $0$ or $\eta$.
However, it turns out there exists a function continuous on the set of matrices with
no singular values equal to $\eta$, which is equal to the mapping
$X \mapsto \left(DF(X)\right)^* Z$, wherever that mapping is defined.
We refer to this function as a continuous extension.
In this section, we find the continuous extension of $X \mapsto \left(DF(X)\right)^* Z$
for all fixed $Z$.

We assume that $X, Z, \Sigma, \zeta, \Gamma, \Delta \in \reals^{m \times n}$,
$U \in \reals^{m \times m}$, $V \in \reals^{n \times n}$, and that
$\Omega_U, \Omega_V, \Omega_\Sigma$
are linear operators from $\reals^{m\times n}$ to $\reals^{m\times n}$.
Without loss of generality, we assume $m \geq n$.
A simple matrix is one without repeated singular values.

\subsection{Gradient for full-rank and simple matrices}

\cite{candes2013svt} gives that for simple and full-rank $X$:
\[
    \left(DF(X)\right)\Delta = U\left(
        (\Omega_{U}\Delta) F(\Sigma) +
        (\Omega_\Sigma\Delta) +
        F(\Sigma)(\Omega_V\Delta)
    \right) V^T
\]
where
\[
    \begin{aligned}
        \left(\Omega_{U}\Delta\right)_{ij} &= \begin{cases}
        0 & \text{if }  i = j \\
        \begin{array}{ll}
            -\frac{1}{\sigma_i^2 - \sigma_j^2} \big(&
            \hspace{-2ex}\sigma_j (U^T \Delta V)_{ij} \\ &+ \sigma_i (U^T \Delta V)_{ji}
        \big)
        \end{array}
            & \text{if } i \neq j \land i \leq n \\
        \frac{1}{\sigma_j} (U^T \Delta V)_{ij} & \text{else } 
    \end{cases}, \\
        \left(\Omega_{V}\Delta\right)_{ij} &= \begin{cases}
        0 & \text{if }  i = j \\
        \begin{array}{ll}
            \frac{1}{\sigma_i^2 - \sigma_j^2} \big(&
            \hspace{-2ex}\sigma_i (U^T \Delta V)_{ij} \\ &+ \sigma_j (U^T \Delta V)_{ji}
        \big)
        \end{array}
   & \text{else } \\
    \end{cases},
    \end{aligned}
\]
and
\[
    \left(\Omega_\Sigma\Delta\right)_{ij} = \begin{cases}
        \mathcal T_\eta'(\sigma_i)(U^T \Delta V)_{ii} & \text{if } i = j \\
        0 & \text{if }  i \neq j \\
    \end{cases}.
\]

In order to find the adjoint of this mapping we begin by constructing a convenient
orthonormal basis of $\reals^{m\times n}$.
We then project the desired quantity $(DF(X))^* Z$ onto the basis vectors.
We can then weight and sum the basis elements to form $(DF(X))^* Z$.

Let $\{E^{ij}\}_{i,j \in [m] \times [n]}$ be the standard basis of $\reals^{m\times n}$,
\ie, $E_{k\ell}^{ij} = 1$ iff $i = k$ and $j = \ell$ and is otherwise $0$.
Let $\Delta^{ij} = u_i v_j^T$.
Critically, $U^T\Delta^{ij} V = E^{ij}$ which will greatly simplify the mappings given above.
For notational simplicity, let $\zeta =U^T Z V$.

Evaluating the projection yields
\[
    \langle (DF(X))^* Z \mid \Delta^{ij} \rangle = 
    \begin{cases}
        \mathcal T_\eta'(\sigma_i) \zeta_{ii} &\text{if } i = j \\
        \frac{\mathcal T_\eta(\sigma_j)}{\sigma_j} \zeta_{ij}   & \text{if } i > n \\
        \begin{array}{l}
        \frac{
            \sigma_i \mathcal T_\eta(\sigma_i) - \sigma_j \mathcal T_\eta(\sigma_j)
        }{\sigma_i^2 - \sigma_j^2}\zeta_{ij} \\
        + \frac{
            \sigma_j \mathcal T_\eta(\sigma_i) - \sigma_i \mathcal T_\eta(\sigma_j)
        }{\sigma_i^2 - \sigma_j^2}\zeta_{ji}
        \end{array}
        &\text{else } \\
    \end{cases}.
\]
This projection is not defined for some basis elements whenever there exists
$i \neq j$ such that $\sigma_i = \sigma_j$ or $i$ such that $\sigma_i = 0$.

\subsection{Extension by continuity to all matrices}
Following \cite{candes2013svt}, we seek to extend the projection of
$(DF(X))^* Z$ by continuity to the situation where there exists $i \neq j$,
such that $\sigma_i = \sigma_j$ or there exists $\sigma_i = 0$.
Note that the projection is only ill-defined for basis elements $\Delta^{ij}$
such that $i \leq n$ and $i \neq j$.
Since simple and full-rank matrices are dense in $\reals^{m\times n}$, we will
consider a sequence of matrices $X^{(k)}$ such that each $X^{(k)}$ is simple and
full-rank and $\lim_{k\to\infty}X^{(k)} = X$.

From \cite{candes2013svt}, we have that for $i \neq j$ such that
$\sigma_i = \sigma_j > 0$,
\[
    \frac{
        \sigma_i^{(k)} \mathcal T_\eta(\sigma_i^{(k)}) - \sigma_j^{(k)} \mathcal T_\eta(\sigma_j)
    }{\left(\sigma_i^{(k)}\right)^2 - \left(\sigma_j^{(k)}\right)^2}\zeta_{ij}\to
    \left(
        \frac{1}{2}\mathcal T_\eta'(\sigma_i) + \frac{1}{2} \frac{\mathcal T_\eta(\sigma_i)}{\sigma_i}
    \right) \zeta_{ij},
\]
and that for $i \neq j$ such that $\sigma_i = \sigma_j = 0$,
\[
    \frac{\sigma_i^{(k)} \mathcal T_\eta(\sigma_i^{(k)}) - \sigma_j^{(k)} \mathcal T_\eta(\sigma_j)}{
        \left(\sigma_i^{(k)}\right)^2 - \left(\sigma_j^{(k)}\right)^2}\zeta_{ij}\to
    \mathcal T_\eta'(0)\zeta_{ij}.
\]
A symmetric version of the argument from \cite{candes2013svt} gives that for $i \neq j$ such that
$\sigma_i = \sigma_j > 0$,
\[
    \frac{\sigma_j^{(k)} \mathcal T_\eta(\sigma_i^{(k)}) - \sigma_i^{(k)} \mathcal T_\eta(\sigma_j)}{
        \left(\sigma_i^{(k)}\right)^2 - \left(\sigma_j^{(k)}\right)^2}\zeta_{ji}\to
    \left(\frac{1}{2}\mathcal T_\eta'(\sigma_i) 
    - \frac{1}{2} \frac{\mathcal T_\eta(\sigma_i)}{\sigma_i} \right) \zeta_{ji},
\]
and for $i \neq j$ such that $\sigma_i = \sigma_j = 0$,
\[
    \frac{\sigma_j^{(k)} \mathcal T_\eta(\sigma_i^{(k)}) - \sigma_i^{(k)} \mathcal T_\eta(\sigma_j)}{
        \left(\sigma_i^{(k)}\right)^2 - \left(\sigma_j^{(k)}\right)^2}\zeta_{ji}\to
    0.
\]
Lastly, note that when $\sigma_j = 0$,
\[
    \lim_{\sigma_j^{(k)} \to 0} \frac{\mathcal T_\eta(\sigma_j^{(k)})}{\sigma_j^{(k)}} = \mathcal T_\eta'(0).
\]
In summary, the continuous extension of $\langle (D F(X))^* Z, \Delta^{ij}$ for
all $X$ is given by
\[
    \Gamma_{ij} = 
\begin{cases}
         \mathcal T_\eta'(\sigma_i) \zeta_{ii} &\text{if } i = j \\
        R(\sigma_j) \zeta_{ij} &\text{if } i > n \\
        Q(\sigma_i, \sigma_j) \zeta_{ij} + T(\sigma_i, \sigma_j) \zeta_{ji}
        &\text{else} \\
    \end{cases}
\]
where
\[
    R(\sigma) = \begin{cases}
        \frac{ \mathcal T_\eta(\sigma)}{\sigma} &\text{if  } \sigma > 0 \\
         \mathcal T_\eta'(\sigma) &\text{if  } \sigma = 0 \\
        \end{cases},
\]
\[
Q(\sigma_i, \sigma_j) = \begin{cases}
        \frac{1}{2}  \mathcal T_\eta'(\sigma_i) 
            + \frac{1}{2} \frac{\mathcal T_\eta(\sigma_i)}{\sigma_i}
        & \text{if } \sigma_i = \sigma_j > 0 \\
         \mathcal T_\eta'(0) &\text{if } \sigma_i = \sigma_j =  0 \\
        \frac{
            \sigma_i  \mathcal T_\eta(\sigma_i) - \sigma_j  \mathcal T_\eta(\sigma_j)
        }{\sigma_i^2 - \sigma_j^2} & \text{else}
    \end{cases},
\]
and
\[
T(\sigma_i, \sigma_j) = \begin{cases}
        \frac{1}{2}  \mathcal T_\eta'(\sigma_i) 
            - \frac{1}{2} \frac{\mathcal T_\eta(\sigma_i)}{\sigma_i}
        & \text{if } \sigma_i = \sigma_j > 0 \\
        0 & \text{if } \sigma_i = \sigma_j =  0 \\
        \frac{
            \sigma_j  \mathcal T_\eta(\sigma_i) - \sigma_i  \mathcal T_\eta(\sigma_j)
        }{\sigma_i^2 - \sigma_j^2} & \text{else}
    \end{cases}.
\]

\subsection{Numerically stable computation}
Constructing $\Delta^{ij}$ in order to evaluate \(\sum_{i=1}^m \sum_{j=1}^n
\Gamma_{ij} \Delta^{ij}\) is numerically unstable in high dimensions.

However, some simple algebra gives that 
\begin{align*}
     \sum_{i=1}^m \sum_{j=1}^n \Gamma_{ij} \Delta^{ij} &=
     U U^T \left(\sum_{i=1}^m \sum_{j=1}^n \Gamma_{ij} \Delta^{ij}\right) V V^T \\ &= 
     U \left(\sum_{i=1}^m \sum_{j=1}^n \Gamma_{ij} U^T \Delta^{ij} V\right) V^T \\ &= 
     U \left(\sum_{i=1}^m \sum_{j=1}^n \Gamma_{ij} E^{ij} \right) V^T = 
     U \Gamma V^T.
\end{align*}
Experimentally, evaluating $U \Gamma V^T$ is numerically stable.

\vfill

\end{document}